\documentclass[final,leqno,onefignum,onetabnum]{siamltex1213}

\usepackage[]{graphicx}

\def\comment#1{}

\newcommand{\bR}{{\bf R}}

\newcommand{\bF}{{\bf F}}
\newcommand{\bY}{{\bf Y}}

\newcommand{\bA}{{\bf A}}

\newcommand{\bG}{{\bf G}}

\newcommand{\bI}{{\bf I}}

\newcommand{\bX}{{\bf X}}
\newcommand{\bE}{{\bf E}}

\newcommand{\bB}{{\bf B}}
\newcommand{\bQ}{{\bf Q}}

\newcommand{\bC}{{\bf C}}

\newcommand{\tS}{\cal{S}}

\newcommand{\tT}{{\cal T}}
\newcommand{\tE}{{\cal E}}

\newcommand{\tM}{{\cal M}}

\title{Non-Orthogonal Tensor Diagonalization % and Block Diagonalization
\thanks{This work was supported by
the Czech Science Foundation through Project No. 14-13713S.}}

\author{Petr Tichavsk\'{y}\thanks{Institute of Information Theory and Automation,
Pod vod\'{a}renskou v\v{e}\v{z}\'{\i} 4, P.O.Box 18, 182 08 Prague 8, Czech Republic. (\email{tichavsk@utia.cas.cz})}}
\author{Petr Tichavsk\'{y}, Anh Huy Phan and Andrzej
Cichocki\thanks{Institute of Information Theory and Automation,
Prague, Czech Republic; Brain Science Institute, RIKEN, Wakoshi,
Japan. A. Cichocki is also with the Systems Research Institute,
Polish Academy of Sciences. }}

\begin{document}
\maketitle
\slugger{mms}{xxxx}{xx}{x}{x--x}%slugger should be set to mms, siap, sicomp, sicon, sidma, sima, simax, sinum, siopt, sisc, or sirev

\begin{abstract}
Tensor diagonalization means transforming a given tensor to an exactly or nearly diagonal form
through multiplying the tensor by non-orthogonal invertible matrices along
selected dimensions of the tensor. It is generalization of approximate joint diagonalization (AJD)
of a set of matrices. In particular, we derive (1) a new algorithm for symmetric AJD, which is called two-sided
symmetric diagonalization of order-three tensor, (2) a similar algorithm for non-symmetric AJD, also called
general two-sided diagonalization of an order-3 tensor, and (3) an algorithm for three-sided diagonalization of order-3 or order-4 tensors.
The latter two algorithms may serve for canonical polyadic (CP) tensor decomposition,
and they can outperform other CP tensor decomposition methods in terms of computational speed
under the restriction that the tensor rank does not exceed the tensor
multilinear rank. Finally, we propose (4) similar algorithms for tensor block diagonalization,
which is related to the tensor block-term decomposition.
\end{abstract}



\pagestyle{myheadings}
\thispagestyle{plain}
\markboth{NON-ORTHOGONAL TENSOR DIAGONALIZATION}{NON-ORTHOGONAL TENSOR DIAGONALIZATION}%{A TOOL FOR BLOCK TENSOR DECOMPOSITIONS}

\section{Introduction}

The approximate joint diagonalization (AJD) of a set of matrices has been recently
recognized to be instrumental in signal processing, mainly
because of its importance in practical signal processing problems
such as source separation, blind beamforming, image denoising,
blind channel identification for multiple-input, multiple-output
(MIMO) telecommunication system, Doppler-shifted echo
extraction in radar, and ICA \cite{SPM}.

Perhaps one of the first such algorithms
is the joint approximate diagonalization of eigenmatrices
(JADE) algorithm proposed in [8]. In this algorithm, the matrices
under consideration are Hermitian and the
considered joint diagonalizer is a unitary
matrix. More recently, generalizations
and/or new decompositions were
found to be of considerable
interest. They concern new
sets of matrices, a nonunitary
joint diagonalizer, and
new decompositions.

The set of given matrices to be diagonalized is a tensor.
The AJD problem can be viewed as a special case of the tensor diagonalization, as we show
later in this paper.

The concept of tensor diagonalization was first introduced by P.
Comon and his co-workers \cite{TDC1,TDCrep}. It works for order-three tensors
of a cubic shape. The tensor
diagonalization in those papers was orthogonal: it
sought orthogonal matrices that would transform the given tensor
in a diagonal one. The method was based on Jacobi rotations.

The tensor diagonalization studied in this paper is
non-orthogonal. We consider two-sided tensor diagonalization of order-three tensors,
which can be symmetric and nonsymmetric, and three-sided diagonalization of order-3 or order-4 tensors.
All algorithms in this paper are based on the same principle.
The main idea is similar to the idea of an AJD algorithm UWEDGE \cite{wasobi},
it can be described in words as ``diagonalize until a further diagonalization is not possible",
but the implementation and performance are different.

In the case of symmetric diagonalization of order-three tensors, we obtain a novel method of AJD.
In the case of nonsymmetric two-sided diagonalization of order-three tensors, we obtain a novel method of
canonical polyadic (CP) tensor decomposition, which follows the idea of SECSI framework for CP decomposition
\cite{SECSI}.

The cases of three-sided diagonalization of order-3 and order-4 tensor represent another method of CP decomposition of order-3 tensor,
and joint approximate diagonalization of several order-3 tensors, respectively.
A generalization to four-sided and more-sided diagonalization of higher-order tensors is straightforward.

The tensor diagonalization methods considered in this paper can be easily modified for block diagonalization.
In many applications, the ordinary diagonalization is not quite appropriate, and
like in independent subspace analysis \cite{ISA}, one seeks rather for
subspaces of columns that represents multidimensional signal components that should be separated or eliminated.
The joint block
diagonalization of the set of these matrices was studied e.g., in \cite{HG,Nion,Lahat,bloky,Shi}.
In the area of tensor decompositions, we speak about the block-term
decomposition, promoted by De Lathauwer and his co-workers
\cite{BTD2,BTD3}. The decomposition means that the given tensor is
rewritten as a sum of several tensors of the same size but a lower
multilinear rank. The block term decomposition was used to propose
a blind DS-CDMA receiver in \cite{Nion-Lath-IEEE}.

In practice, initializing a BTD without getting captured in false local minima of the criterion
function is a very challenging
problem. Another difficulty is that the appropriate block sizes
might not be known in advance. In some cases we have empirically found that the
tensor diagonalization can be used to carry out a suitable
block-term decomposition, i.e., find appropriate block sizes, provided there is no or little noise.
This has been already observed in \cite{LVAICA}.

There are a few related conference publications on the topic. The original version of the paper considered
tensor diagonalization through generalized Jacobi (Givens) rotations \cite{arxiv}. An algorithm for
two-sided diagonalization of order-3 tensor was proposed in a conference paper \cite{Shen}.
An algorithm for three-sided diagonalization of order-3 tensor was
proposed in \cite{Maurandi}.
This paper presents a different approach to the same problem.

The paper is organized as follows: Section 2 presents the basic
principles of tensor diagonalization and shows its connection to
CP decomposition. In Section 3, iterative algorithms are proposed to perform the
three-sided and two-sided symmetric and nonsymmetric diagonalization, either in the real or complex domain. In section 4,
new algorithms for joint block diagonalization and block-term decomposition are developed.
Section 5 presents some numerical examples, and Section 6 concludes the paper.

\section{Tensor diagonalization principle}

The main idea of the tensor diagonalization is to find so-called
demixing matrices that transform the given tensor to another
tensor that is diagonally dominant.
In the AJD, we are given a set of matrices $\bR_m$, $m=1,\ldots,M$
and we seek for so-called de-mixing matrix $\bA$ such that
$\bA\bR_m\bA^H$, $m=1,\ldots,M$ are all diagonally dominant. It means that
the off-diagonal elements are significantly smaller in magnitude than
the diagonal elements. Here, $^H$ denotes the Hermitian transpose.

There are several measures of success and several algorithms that accomplish the diagonalization,
see \cite{SPM} for a review. Some of them can be modified to provide approximate joint block diagonalization.

A modification of the problem is the nonsymmetric AJD. Here we assume again, that the given matrices
$\bR_m$, $m=1,\ldots,M$ are square, and we seek for invertible matrices $\bA$, $\bB$ such that
$\bA\bR_m\bB^T$, $m=1,\ldots,M$ are all diagonally dominant. The symbol $^T$ stands for the matrix transpose.
A tensor formulation of the same problem can be following:

Let $\tT$ be a tensor of size $n\times n\times m$ composed of the slices $\{\bR_m\}$, $m=1,\ldots,M$.
The outcome of the diagonalization is the tensor
$$
\tE=\tT\times_1 \bA\times_2 \bB
$$
where $\times_i$ denotes the tensor-matrix multiplication along the dimension $i$, $i=1,2$.
%, and $^*$ is the
%complex conjugate. The complex conjugation here is immaterial, but we include it here to make the
%problem equivalent to the non-symmetric AJD.
A successful diagonalization means that $\|\mbox{off}_2(\tE)\|_F$ is small compared to diagonal elements of $\tE$,
% can be measured by the ratio $\|\mbox{off}_2(\tE)\|_F^2/\|\tE\|_F^2$,
where $\|\cdot\|_F$ is the Frobenius norm,
and  $\mbox{off}_2(\tE)$ is the operator that nullifies all elements of the input tensor
except the diagonals of all frontal slices of the tensor. In other words, if
$\tE$ has elements $\tE_{ijk}$, then $\mbox{off}_2(\tE)$ has elements $(1-\delta_{ij})\tE_{ijk}$,
where $\delta_{ij}$ is the Kronecker delta.

Similarly, three-sided diagonalization of an order-4 tensor $\tT$ of the
size $N\times N\times N\times M$ with elements
$t_{ijkm}$, $i,j,k=1,\ldots,N$, $m=1,\ldots,M$ consists in finding three matrices $\bA$, $\bB$ and $\bC$
of size $N\times N$ such that
\begin{equation}
\tE=\tT\times_1 \bA\times_2 \bB \times_3 \bC \label{TD}
\end{equation}
is nearly {\em spatially} diagonal in the sense
$$
\|\mbox{off}_3(\tE)\|_F^2\ll \|\tE\|_F^2
$$
where the operator $\mbox{off}_3$ nullifies all elements of $\tE$ that do not lie on the spatial diagonal of the tensor.
To be exact, the operator $\mbox{off}_3$ acts on a tensor $\tE$ with elements $\tE_{ijkm}$ so that
$\mbox{off}_3(\tE)$ has elements $(1-\delta_{ij}\delta_{jk})\tE_{ijkm}$.
In the special case $M=1$, $\tT$ and $\tE$ are order-3 tensors because of having only three variable indices.
The diagonalization is illustrated in Fig. 1 for the case $M=1$. The
multiplication in (\ref{TD}) can be written as
\begin{equation}
e_{ijkm} =\sum_{\alpha,\beta,\gamma=1}^N t_{\alpha\beta\gamma m}
a_{i\alpha} b_{j\beta} c_{k\gamma}
\end{equation}
where $e_{ijkm}$, $t_{\alpha\beta\gamma m}$, $a_{i\alpha}$,
$b_{j\beta}$ and $c_{k\gamma}$ are elements of tensors $\tE$,
$\tT$ and matrices $\bA$, $\bB$ and $\bC$, respectively.

%For easy reference, consider the off-diagonality operator $\mbox{off}(\cdot)$ which transforms
%a tensor $\tE$ of the size $N\times N\times N\times M$ with elements $e_{ijkm}$
%into a tensor of the same size, having the elements $(1-\delta_{ij}\delta_{jk})e_{ijkm}$ for
%$i,j,k=1,\ldots,N$ and $m=1,\ldots,M$, where $\delta_{ij}$ is the Kronecker delta.
%In other words, the operator nullifies the elements $e_{iiim}$, $i=1,\ldots,N$ and $m=1,\ldots,M$,
%and keeps the remaining tensor elements unchanged.
%
%The approximate diagonalization means seeking for matrices $\bA$,
%$\bB$ and $\bC$ of a size $N\times N$ such that the product
%\begin{equation}\label{tedia}
%\tE =\tT \times_1 \bA \times_2 \bB \times_3 \bC
%\end{equation}
%is approximately spatially diagonal in the sense
%$\mbox{off}(\tE)\approx 0$.

The success of the diagonalization can be defined in different ways. The algorithms proposed in this paper
are based on the principle ``diagonalize until a further diagonalization is not possible".
%This principle implies a measure of success which is unique.
Let us explain the principle on the three side diagonalization.

The condition that the resulting tensor $\tE$ ``cannot be
diagonalized more" means that
\begin{equation}\label{tedia1}
\|\mbox{off}_3(\tE^\prime)\|_F \geq \|\mbox{off}_3(\tE)\|_F
\end{equation}
for all
\begin{equation}\label{tedia2}
\tE^\prime =\tE \times_1 \widetilde\bA \times_2 \widetilde\bB \times_3 \widetilde\bC
\end{equation}
where diagonals of $\widetilde\bA$, $\widetilde\bB$, $\widetilde\bC$ are filled with 1's,
symbolically
\begin{equation}\label{tedia3}
\mbox{diag}(\widetilde\bA)=\mbox{diag}(\widetilde\bB)=\mbox{diag}(\widetilde\bC)=\mbox{diag}({\bf I})=(1,\ldots,1)^T~.
\end{equation}
The objective function to be minimized is the norm of the gradient of $\|\mbox{off}(\tE^\prime)\|_F$
with respect to vector of off-diagonal elements of $\widetilde\bA$, $\widetilde\bB$, and $\widetilde\bC$
at the point $\widetilde\bA=\widetilde\bB=\widetilde\bC=\bI$~.
Ideally, the norm of the gradient should be zero, and the corresponding Hessian matrix should be positive definite.
\begin{figure}
\centering
\includegraphics[width=0.7\linewidth, trim = 0.0cm .0cm 0cm 0cm,clip=false]{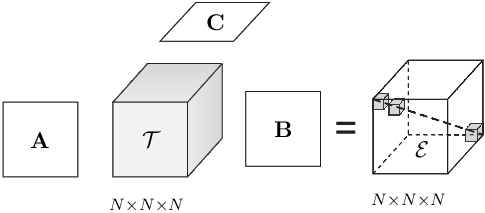}

\vspace{1ex}
\caption{Three-sided tensor diagonalization
transforms a tensor $\tT$ to a diagonally dominant tensor $\tE$
using demixing matrices $\bA$, $\bB$ and $\bC$: $\tT  \times_1 \bA
\, \times_2 \, \bB \, \times_3 \, \bC = \tE$.}
\end{figure}
%Note that the similar idea was used in the algorithm UWEDGE and WEDGE for joint approximate diagonalization
%of a set of matrices \cite{uwedge}.

The tensor diagonalization is not unique. Like in the CP
decomposition, there is a permutation ambiguity, meaning that the
order of rows in $\bA$ and accordingly in $\bB$ and $\bC$ can be
arbitrary. Moreover, there is also a scale ambiguity, if $\mbox{off}_3(\tE)={\bf 0}$,
i.e. when the tensor admits an exact CP decomposition, and all factor matrices are invertible.
In that case, we will show that the tensor diagonalization is essentially unique, and its outcome is
equivalent to the CP decomposition.
In other cases, the tensor diagonalization is not unique. The diagonalization might have several (perhaps
infinitely many) possible outcomes, but any of them characterizes the tensor in a sense, and may reveal
hidden block structure in the tensor. What we typically observe that the output core tensor $\tE$ contains
many nulls (entries with negligible magnitudes) and is sparse in this sense. %The purpose is the exploratory analysis of the tensor.

\section{Algorithm TEDIA}

The tensor diagonalization can be achieved by a cyclic
application of elementary rotations for all pairs of distinct
indices $i,j=1,\ldots,N$, as it was shown in earlier versions of this paper \cite{arxiv}.
Here, however, we present an easier way to achieve the goal.
The proposed method is a gradient method with an enhanced line search, similar to \cite{ELS}.
We explain it in the case of three-sided diagonalization first.

\subsection{Three-Sided Diagonalization}

Assume that $\tE$ is a partially diagonalized tensor obtained
during the optimization process.
Let $\bG_A$ be the gradient of the function
$\|\mbox{off}(\tE\times_1 \bA)\|_F^2$ with respect to $\bA$ at $\bA=\bI$.
Similarly, let $\bG_B$ and $\bG_C$ be gradients o
$\|\mbox{off}(\tE\times_2 \bB)\|_F^2$ with respect to $\bB$ at $\bB=\bI$, and of
$\|\mbox{off}(\tE\times_3 \bC)\|_F^2$ with respect to $\bC$ at $\bC=\bI$, respectively.
The diagonal elements of $\bG_A$, $\bG_B$, and $\bG_C$ are set to zero, because
the diagonals of $\bA$, $\bB$ and $\bC$ are fixed.

It can be easily found (see Appendix A) that
\begin{eqnarray}
\bG_A & = & \mbox{off}[(\mbox{off}_3(\tE))_{(1)}\tE_{(1)}^T]\nonumber \\
\bG_B & = & \mbox{off}[(\mbox{off}_3(\tE))_{(2)}\tE_{(2)}^T]\label{gradient}\\
\bG_C & = & \mbox{off}[(\mbox{off}_3(\tE))_{(3)}\tE_{(3)}^T]\nonumber
\end{eqnarray}
where $\tE_{(i)}$ and $(\mbox{off}_3(\tE))_{(i)}$ are the mode-$i$ matricizations of
$\tE$ and $\mbox{off}_3(\tE)$, respectively.

Once $\bG_A$, $\bG_B$ and $\bG_C$ are computed, we seek a scalar step size $t$
which minimizes the following polynomial of degree 6,
\begin{equation}
\varphi(t)=\|\mbox{off}_3(\tE\times_1 (\bI+t\bG_A)\times_2 (\bI+t\bG_B)\times_3 (\bI+t\bG_C)\|_F^2~.\label{polynom}
\end{equation}
%Note that the  minimum of $\varphi(t)$ can be found among roots of
%the first derivative of the polynomial.
Let $t_m$ be the minimizer of $\varphi(t)$.
Then, the next iteration of $\tE$ is obtained as
\begin{equation}
\tE \leftarrow % \tE^\prime =
\tE\times_1 (\bI+t_m\bG_A)\times_2 (\bI+t_m\bG_B)\times_3 (\bI+t_m\bG_C)~.\label{update1}
\end{equation}
The estimated demixing matrices are updated as
\begin{eqnarray}
\bA &\leftarrow& %\bA^\prime =
(\bI+t_m\bG_A)\bA\nonumber\\
\bB &\leftarrow& %\bB^\prime =
(\bI+t_m\bG_B)\bB \label{update2}\\
\bC &\leftarrow& %\bC^\prime =
(\bI+t_m\bG_C)\bC\nonumber
\end{eqnarray}
The algorithm is summarized in Table 1. In the table, the notation $<\tE_i,\tE_j>$ denotes a scalar product of tensors
$\tE_i$, $\tE_j$, i.e., sum of entries of the elementwise product of these tensors.
%\end{document}
\begin{table}[t]
\caption{Algorithm TEDIA for three-sided diagonalization}
{\bf Input:}  Tensor $\tT$ of size $N\times N\times N\times M$, stopping constant $\varepsilon$\\
{\bf Output:} Mixing matrices $\tilde{\bA}$, $\tilde{\bB}$, $\tilde{\bC}$; Core tensor $\tE$\\
Initialize: $\tE:=\tT$, $\tilde{\bA}=\tilde{\bB}=\tilde{\bC}=\bI$ \\
Repeat \\
\rule{0mm}{0mm}\quad 1) Compute the gradients $\bG_A,\bG_B,\bG_C$ in (\ref{gradient})\\
\rule{0mm}{0mm}\quad 2) Compute coefficients of the polynomial (\ref{polynom})
\begin{eqnarray*}
\varphi(t)&=& c_0+c_1 t+c_2t^2+c_3 t^3+c_4 t^4 +c_5 t^5+c_6 t^6\\
\tE_1 &=& \mbox{off}_3(\tE\times_1 \bG_A +\tE \times_2 \bG_B +\tE \times_3 \bG_C) \\
\tE_2 &=& \mbox{off}_3(\tE\times_1 \bG_A \times_2 \bG_B+\tE \times_2 \bG_B\times_3 \bG_C +\tE \times_1 \bG_A\times_3 \bG_C) \\
\tE_3 &=& \mbox{off}_3(\tE\times_1 \bG_A \times_2 \bG_B \times_3 \bG_C) \\
c_0 &=& <\mbox{off}_3(\tE),\tE>=\|\mbox{off}_3(\tE)\|_F^2\\
c_1 &=&  2<\tE_1,\tE>\\
c_2 &=& \|\mbox{off}_3(\tE_1)\|_F^2+2<\tE_2\tE>\\
c_3 &=& 2<\tE_1,\tE_2>+2<\tE_3,\tE>\\
c_4 &=& \|\mbox{off}_3(\tE_2)\|_F^2+2<\tE_1,\tE_3>\\
c_5 &=& 2<\tE_2,\tE_3>\\
c_6 &=& \|\mbox{off}_3(\tE_3)\|_F^2
\end{eqnarray*}
%\rule{0mm}{0mm}\qquad where $<\tE_i,\tE_j>$ is the sum of elements of $\mbox{off}_3(\tE_i *\tE_j)$.\\
\rule{0mm}{0mm}\quad 3) Find the root $t_m$ of $\varphi^\prime(t)=c_1+2c_2t+3c_3 t^2+4c_4 t^3+5c_5 t^4+6c_6 t^5$ minimizing $\varphi(t)$\\
\rule{0mm}{0mm}\quad 4) Update $\tE$, $\bA, \bB, \bC$ as in (\ref{update1}) and (\ref{update2})~.\\
\rule{0mm}{0mm}Until $\|\bG_A\|_F+\|\bG_B\|_F+\|\bG_C\|_F<\varepsilon$
\end{table}

The computational complexity of the TEDIA algorithm is $O(N^4M)$ operations per iteration. In the case $M=1$, the complexity per iteration is
roughly the same as the complexity of one iteration of the Alternating Least Squares (ALS) algorithm with the Enhanced Line Search (ELS).
The convergence of TEDIA appears to be smoother than that of ALS or ALS/ELS, namely in difficult scenarios, as we show in the simulation section.

\subsection{Two-Sided Diagonalization}

A modification of TEDIA to the two-sided diagonalization is straightforward. The difference is that
\begin{eqnarray}
\bG_A & = & \mbox{off}[(\mbox{off}_2(\tE))_{(1)}\tE_{(1)}^T]\nonumber \\
\bG_B & = & \mbox{off}[(\mbox{off}_2(\tE))_{(2)}\tE_{(2)}^T]\label{gradient2}
\end{eqnarray}
and the polynomial $\varphi(t)$ would be of degree 4,
\begin{equation}
\varphi(t)=\|\mbox{off}_2(\tE\times_1 (\bI+t\bG_A)\times_2 (\bI+t\bG_B)\|_F^2~.\label{polynomC}
\end{equation}
In the case of symmetric two-sided diagonalization, there is only one de-mixing matrix $\bA$,
the corresponding gradient matrix is
\begin{eqnarray*}
\bG_A &=& \mbox{off}[(\mbox{off}_2(\tE))_{(1)}\tE_{(1)}^T]
\end{eqnarray*}
the polynomial $\varphi(t)$ would be of degree 4 again,
\begin{equation}
\varphi(t)=\|\mbox{off}_2(\tE\times_1 (\bI+t\bG_A)\times_2 (\bI+t\bG_A)\|_F^2~.\label{polynomD}
\end{equation}
% and to symmetric two-sided diagonalization is straightforward.

\section{Application in CP tensor decompositions}

A natural utilization of TEDIA is in CP tensor decomposition. It was shown in \cite{SECSI} and \cite{icassp16} (so-called SECSI framework) that two-sided
tensor diagonalization can be applied in CP tensor decomposition. Similarly, the three-sided diagonalization can be used for this purpose.
Let us discuss this issue in more details.

First of all, the tensor might have a different shape than $N\times N\times N$ or $N\times N\times N\times M$. A direct tensor diagonalization makes no sense because the diagonalizing matrices must be square and invertible. The Tucker compression is advised in this case. The Tucker compression consists in finding orthogonal matrices $\bQ_1,\bQ_2,\bQ_3$
such that
$$
\tT\approx \tT_C \times_1\bQ_1\times_2\bQ_2\times_3\bQ_3
$$
where $\tT_C$ is the compressed tensor of the required shape, and
$\times_i$ denotes a mode-$i$ tensor-matrix multiplication. The
Tucker compression can be achieved by the HOOI algorithm
\cite{HOOI}, see \cite{HOOI3} for more literature on this topic. %, or \cite{FastTucker}. 
Note that TEDIA can serve as a tool for the Tucker compression as well,
it only would have to be a block version of it, see the next section. However, its performance in the compression
appears to be not that good as the performance of the methods mentioned above.

Also, note that higher-order tensors can also be decomposed through CP decomposition of order-three tensors through the tensor
re-shaping \cite{reshape}. Therefore we shall concentrate on CP decomposition of cubic shaped tensor.

Another remark is on the computational accuracy. The tensor diagonalization is not a statistically efficient procedure of the CP decomposition,
because it is not equivalent to the maximum likelihood estimate. The ultimate performance of TEDIA-based procedure can be achieved
of the TEDIA outcome is taken as input for another technique, which maximizes the likelihood function.
Although the Tucker compression pre-processing is widely used,
it is probably not fully information preserving, and its application prior CP decomposition may result in a
certain loss in accuracy, also. Therefore the maximum likelihood (least squares) estimation should be performed on the original
(not compressed) data.

The theoretical CP decomposition might involve rank-deficient factor matrices. In that case, the optimum de-mixing matrices would not be invertible. This is in conflict with the tensor diagonalization, which always produces invertible demixing matrices.
TEDIA therefore might not be useful in such cases and block TEDIA or block term decomposition would be more appropriate \cite{Stegeman}.

Assume that a tensor $\tT$ is diagonalized by three matrices $\bA$, $\bB$ and $\bC$ such that
the product
$\tE=\tT\times_1\bA\times_2\bB\times_3\bC$ cannot be diagonalized more in the sense of section 2.
Assume the tensor is of size $N\times N\times N$ and the rank is $R\leq N$.
It may occur that the core tensor $\tE$ has only at most $R$ significant nonzero elements, while the magnitude of the other elements
is negligible. Zeroing other than the $R$ significant elements we get a rank-$R$ approximation of the core tensor, which
implies a rank-$R$ approximation of the original tensor.
If the significant elements lie on the main spatial diagonal of the tensor, we have got ordinary diagonalization
and CP decomposition with regular factor matrices. However, if the multilinear rank of the tensor is not $(N,N,N)$, this is not possible,
and not all significant elements lie on the diagonal. In that case, some of the factor matrices in the CP approximation of the tensor
would be rank deficient.

\section{Block diagonalization}

The tensor block diagonalization is a natural generalization of the diagonalization considered in the previous sections.
It can have the form of symmetric or nonsymmetric two-side block diagonalization, or three side diagonalization, as it is
illustrated in Figure 2.
\begin{figure}
\centering
\includegraphics[width=0.7\linewidth, trim = 0.0cm .6cm 0cm
0cm,clip=false]{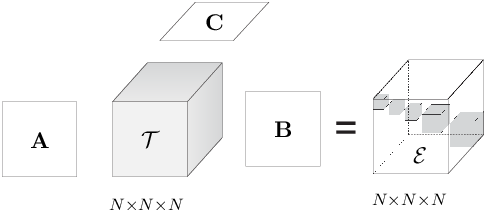}
\vspace{3ex}

\caption{Tensor block diagonalization
transforms a tensor $\tT$ to a block diagonal tensor $\tE$ by
factor matrices $\bA$, $\bB$ and $\bC$: $\tT  \times_1 \bA \,
\times_2 \, \bB \, \times_3 \, \bC =  \tE$.}
\end{figure}
In this section, we assume that the block structure is known.
The block structure of $\tE$ can be represented by a tensor $\tM$ of the same size as $\tE$, which contains zeros
in the place of the desired blocks, and ones elsewhere. In the special case when there is only a single block,
we receive a novel method of the Tucker compression.

We can consider the operator $\mbox{boff}_\tM$ which nullifies all off-block elements of the input tensor,
\begin{equation}
\mbox{boff}_\tM(\tE)=\tE\star\tM
\end{equation}
where $\star$ stands for the elementwise product. Now, we apply the principle
``diagonalize until a further diagonalization is not possible" to get a block diagonalization procedure like in Section 4.
In the case of the three-sided block diagonalization, we define
$\bG_A$, $\bG_B$ and $\bG_C$ as the gradient of the function
$\|\mbox{boff}_\tM(\tE\times_1 \bA)\|_F^2$ with respect to $\bA$ at $\bA=\bI$, gradient of
$\|\mbox{boff}_\tM(\tE\times_2 \bB)\|_F^2$ with respect to $\bB$ at $\bB=\bI$, and gradient of
$\|\mbox{boff}_\tM(\tE\times_3 \bC)\|_F^2$ with respect to $\bC$ at $\bC=\bI$, respectively.
The diagonal elements of $\bG_A$, $\bG_B$, and $\bG_C$ are set to zero, because
the diagonals of $\bA$, $\bB$ and $\bC$ are fixed.
The result is
\begin{eqnarray}
\bG_A & = & \mbox{off}[(\tM\star\tE)_{(1)}\tE_{(1)}^T]\nonumber \\
\bG_B & = & \mbox{off}[(\tM\star\tE)_{(2)}\tE_{(2)}^T]\label{gradientB}\\
\bG_C & = & \mbox{off}[(\tM\star\tE)_{(3)}\tE_{(3)}^T]\nonumber
\end{eqnarray}
Finally, we find the optimum step-size $t_m$ by minimizing the function
\begin{equation}
\varphi(t)=\|\mbox{boff}_\tM(\tE\times_1 (\bI+t\bG_A)\times_2 (\bI+t\bG_B)\times_3 (\bI+t\bG_C)\|_F^2~.\label{polynomB}
\end{equation}
Then, the core tensor $\tE$ and the de-mixing matrices $\bA,\bB$ and $\bC$ are updated as in (\ref{update1}) and (\ref{update2}), respectively.

\section{Blind Block Diagonalization}

By blind block diagonalization we understand block diagonalization
without knowing the block structure in advance.
In \cite{LVAICA} it was shown that an ordinary approximate joint
diagonalization algorithm for matrices can be used to obtain a
joint block diagonalization of these matrices. It appears that
similar link exists between the tensor diagonalization and tensor
block-term decomposition. It is possible that a given tensor may
not be fully diagonalizable but may still be
block diagonalizable, as it is shown schematically in Fig. 3. We
can assume that the diagonal blocks cannot be diagonalized
further, because their tensor ranks exceed their dimensions, and
that each of the blocks separately obeys the zero gradient condition
$\bG_A=\bG_B=\bG_C={\bf 0}$. It is straightforward to prove
that a compound block diagonal tensor obeys the condition as well.

Note that the order of rows in matrices $\bA$, $\bB$ and $\bC$
might be arbitrary, all providing equivalent diagonalizations.
In other words, if $\pi(\cdot)$ is an
arbitrary permutation of $(1,\ldots,N)$, then
\begin{equation}\label{tediaprime}
\tE^\prime =\tT \times_1 \bA^\prime \times_2 \bB^\prime \times_3
\bC^\prime
\end{equation}
where $\tE^\prime$, $\bA^\prime$, $\bB^\prime$ and $\bC^\prime$
have elements $e^\prime_{ijkm}=e_{\pi(i),\pi(j),\pi(k),m}$,
$a^\prime_{i,\alpha}=a_{\pi(i),\alpha}$,
$b^\prime_{j,\beta}=b_{\pi(j),\beta}$, and
$c^\prime_{k\gamma}=c_{\pi(k),\gamma}$, respectively, for
$i,j,k,\alpha,\beta,\gamma=1,\ldots,N$, is an equivalent
diagonalization.
It follows that if the block structure exists in the original
tensor, the block structure may not be apparent after a mixing and
demixing (diagonalization) unless a suitable permutation $\pi$ has
been found. The permutation should be the same in all modes in
order to guarantee that all diagonal elements of the original
tensor will appear on the diagonal of the permuted tensor.

The degree of diagonality or block diagonality of a tensor $\tE$
can be judged via the matrix $\bF$ of size $N\times N$, having
elements
\begin{equation}
f_{ij}=\sum_{m=1}^M\sum_{k=1}^N |e_{kijm}|+|e_{ikjm}|+|e_{ijkm}|\label{FF}
\end{equation}
The tensor $\tE$ is said to be diagonal (block diagonal) if and
only if $\bF$ is diagonal (block diagonal). If $\bF$ is
block diagonal, its $(i,j)-$th element $f_{ij}$ is zero if $i,j$
belong to different blocks, and it might be strictly positive if
they belong to the same block. The order of the columns in the
factor matrices is random in general, however. Then, $f_{ij}$ or
its symmetrized version $f_{ij}+f_{ji}$ may be considered as a
measure of similarity between columns $i$ and $j$ in the mixing
matrices, or an indicator of ``probability" that they belong to
the same block.

Such permutation can be found, e.g., using the well-known reverse
Cuthill-McKee algorithm (RCM)\cite{symrcm}, implemented
in Matlab$^{\rm TM}$ as function {\em symrcm}.
The RCM algorithm, applied to the matrix $\bF$, reveals an ordering of the columns
and rows such that the reordered matrix is block diagonal, if such ordering exists.

In the noisy case, when the blocks of the core tensors are fuzzy, %see Figure 8b,
we have better experiences with standard clustering methods, such
as {\em hierarchical clustering with the average-linking policy}
\cite{clustering}, which take $\bF$ for a similarity matrix. In
short, the algorithm begins with a trivial clustering which
consists of $N$ singletons, and in each subsequent step it merges
those two clusters that have the maximum average similarity
between their members. The algorithm is summarized in Table 2.
\begin{table}[t]
\caption{Clustering of components in MATLAB notations}
{\bf Input:}  Similarity matrix $\bF$ of size $N\times N$ (destroyed in the procedure)\\
{\bf Output:} Permutation $J$ of indices $1,\ldots,N$ such that $\bF(J,J)$
is approximately block diagonal.
$$
\begin{array}{ll}
 \bF(1:N+1:N^2)=0 & \%\,\, \mbox{nullify diagonal of $\bF$, i.e., $\bF(i,i)=0$ for $i=1,\ldots,N$} \\
 S=repmat((1:N)',1,N); &  \%\,\, \mbox{auxiliary array of size $N\times N$, i.e., $S(i,j)=i$} \\%, the $n-$th row is $(n,n,\ldots,n)$
 L=ones(N,1); & \%\,\, \mbox{auxiliary array of the clusters lengths}\\
\mbox{For}\,\, i=N:-1:2 &\%\,\, \mbox{in the $i-th$ step two of $i$ clusters are merged.}\end{array}$$\vspace{-5mm}

\rule{0mm}{0mm}\qquad  $   [m1,\sim]=\max(\bF(1:i,1:i));$\\
\rule{0mm}{0mm}\qquad   $  [m2,j]=\max(m1);$\hfill  \%\,\, $(j,k)$\,\,\mbox{be the clusters with the highest similarity}\\
\rule{0mm}{0mm}\qquad   $  [m3,k]=\max(\bF(:,j));$\\
\rule{0mm}{0mm}\qquad      \mbox{if}\,\, $j>k,\, aux=j;\, j=k;\, k=aux;\,\mbox{end}\hfill$ \%\,\,\mbox{to make sure that $k>j$}\\
\rule{0mm}{0mm}\qquad  $ Lnew=L(j)+L(k);$ \hfill \%\,\, \mbox{length of the new cluster, union of $j,k$}\\
\rule{0mm}{0mm}\qquad  $ Snew=[S(j,1:L(j)),\, S(k,1:L(k)),\, zeros(1,N-Lnew)];$\\
\rule{0mm}{0mm}\qquad    \hfill \%\,\, \ldots\mbox{indices belonging to the new cluster}\\
\rule{0mm}{0mm}\qquad   $ ind=[1:j-1,\,j+1:k-1,\, k+1:i];$\hfill \%\,\,\mbox{indices of the other clusters}\\
\rule{0mm}{0mm}\qquad  $ Fnew=(L(j)*\bF(j,ind)+L(k)*\bF(k,ind))/Lnew;$\\
\rule{0mm}{0mm}\qquad\hfill \%\,\,\mbox{similarities between the new cluster and the other clusters}\\
\rule{0mm}{0mm}\qquad   $ \bF=[0,Fnew;\, Fnew', \bF(ind,ind)];$\hfill \%\,\,\mbox{update of the similarity matrix}\\
\rule{0mm}{0mm}\qquad    $ S=[Snew;\, S(ind,:)];$\hfill \%\,\,\mbox{update indices in the clusters}\\
\rule{0mm}{0mm}\qquad   $ L=[Lnew;\, L(ind)];$\hfill \%\,\,\mbox{update the cluster lengths} \\
\mbox{End}\\
$J=S(1,:);$\\
\mbox{End}
\end{table}

Note that even if the desired block structure of the core tensor is known in advance,
it might be useful to apply the blind diagonalization + clustering as a pre-processing step
for the ordinary (non-blind) block diagonalization, because it may reduce the number of
iterations of the latter algorithm needed to achieve convergence.

\section{Simulations}

\subsection{Example 1: CP decomposition}

In this example, we apply the three-sided tensor diagonalization and other methods to decompose a cubic tensor
of size $20\times 20\times 20$ of rank 20 with collinearity in two and three modes, respectively. The decomposition is hard for all existing methods.
First, we generated three orthogonal matrices of the size $20\times 20$, denoted $\bA_0$, $\bB_0$ and $\bC_0$.
We divided each of them to four blocks of size $20\times 5$, i.e. $\bA_0=[\bA_{01},\bA_{02},\bA_{03},\bA_{04}]$.
Then, the factor matrix $\bA$ was built of four blocks $\bA=[\bA_1,\bA_2,\bA_3,\bA_4]$, where
$\bA_k=c\,\bA_{0k}(:,1){\bf 1}_{1\times 5}+\sqrt{1-c^2}\bA_{0k}$ for $k=1,2,3,4$, $c$ is a free parameter,
$\bA_{0k}(:,1)$ is the first column of $\bA_{0k}$,
and ${\bf 1}_{1\times 5}$ is a row vector of 1's of the size $1\times 5$. We set $c=0.99$.
Thanks to this definition, each of the blocks
$\bA_k$ contained five nearly colinear columns.
Similarly, $\bB$ and $\bC$ were composed of four blocks of nearly co-linear columns obtained using corresponding blocks of $\bB_0$, and $\bC_0$.
Finally, we added an i.i.d. Gaussian noise to each tensor element so that a chosen signal-to-noise ratio (SNR) is attained.
This setting was also used in \cite{PALS}.

We study the performance of six CP decomposition methods. First of all, it is the Direct Tri-Linear Decomposition \cite{DTLD}, which is based on generalized eigendecomposition of a matrix pair. We consider two variants of the method: either we take the
generalized eigendecomposition of the first two frontal slices of the tensor, or generalized eigendecomposition of the two frontal slices of the tensor compressed to the size $N\times N\times 2$. Among the two DTDL results, we consider the one with a lower
fitting error (i.e. Frobenius norm of the tensor and its rank-$N$ CP approximation).
The latter variant is advocated in \cite{DTLD}, the former variant is known to provide an exact solution if there is no noise.
The better of the two variants is referred to as DTLD.
We use this method to initialize all other CP decomposition methods. This algorithm is selected for not giving a large advantage to TEDIA
algorithms compared to the traditional ones. TEDIA is much less sensitive to a wrong initialization because it needs only
a larger number of iterations to achieve convergence if it is initialized wrongly (say randomly). On the other hand, the traditional methods converge to some
local minimum of the quadratic cost function only, and increasing the number of the iterations may not help. In general, DTLD is a fast algorithm.

The second method in the study is the traditional Alternating Least Squares (ALS)
method with 2000 iterations. The third method is the ALS with the Exact Line Search (ELS) method \cite{ELS}, again with 2000 iterations. Fourth, it is the Levenberg-Marquardt (LM) method \cite{SIAM}, where we used 30 iterations. Each iteration is computationally more complex than the other methods, and therefore we have chosen 30 iterations to keep the computational time comparable to the other algorithms. We note that if the number of iterations is increased, say to 2000, the performance of LM would improve, and the method would outperform all other algorithms.

The fifth  method is the three-sided tensor diagonalization method described in this paper, and the sixth method is the
two-sided diagonalization, used like in the SECSI framework \cite{SECSI}. The last two methods stopped after 2000 iterations.
The inverses of the estimated demixing matrices are taken as estimates of the factor matrices.
The two-sided diagonalization produces only estimates of two factor matrices: in this case, the third factor matrix is computed by the least squares fitting as in the ALS method.
Note that each run of the algorithms took 0.07, 1.67, 4.94, 2.48, 3.11 and 1.42 seconds for DTLD, ALS, ALS-ELS, LM, TEDIA 3 and TEDIA 2, respectively.

\begin{figure}
            \centering
            \includegraphics[width=0.48\linewidth]{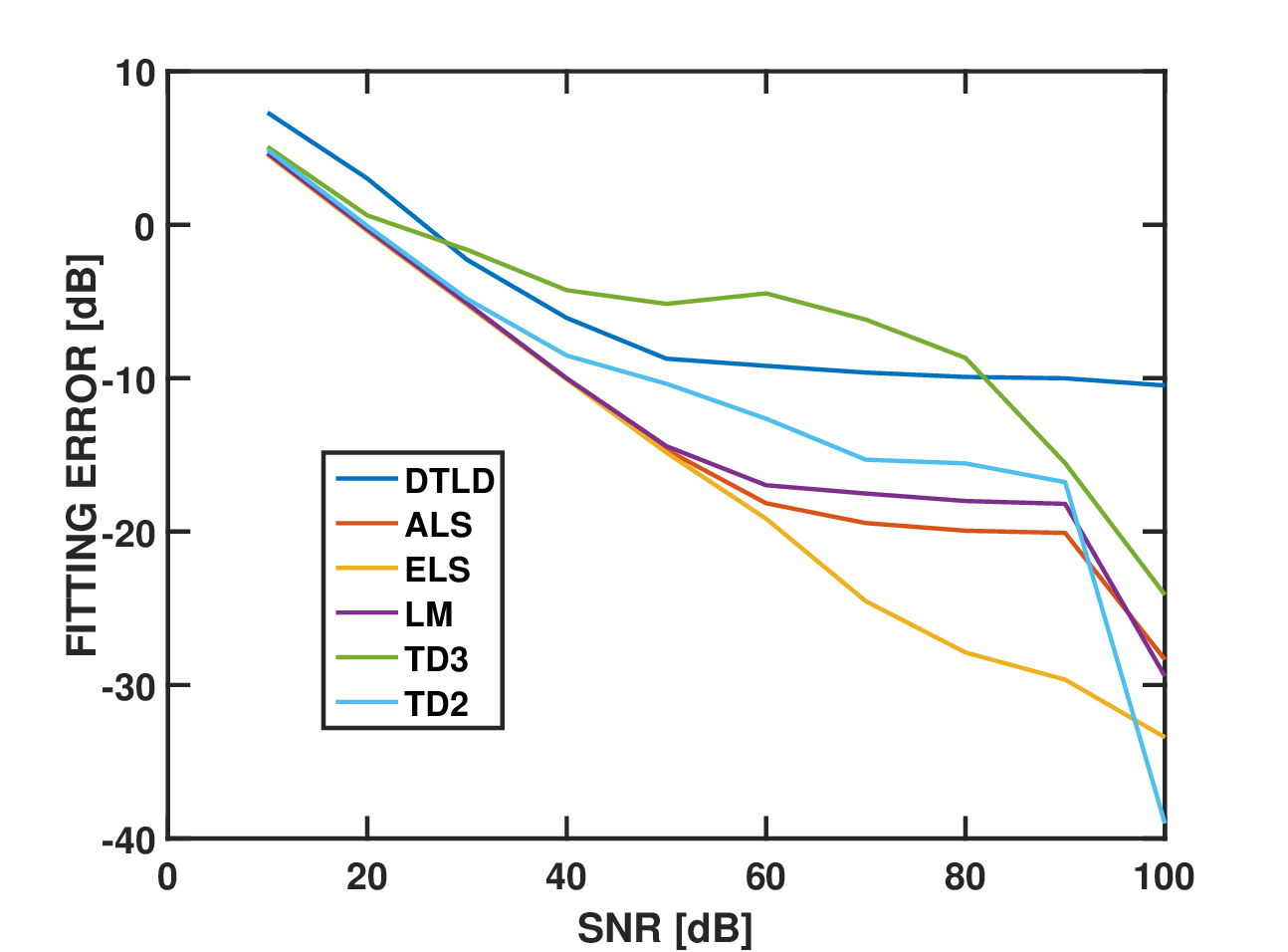}\quad%{bcd_I12R3_msae}
            \includegraphics[width=0.48\linewidth]{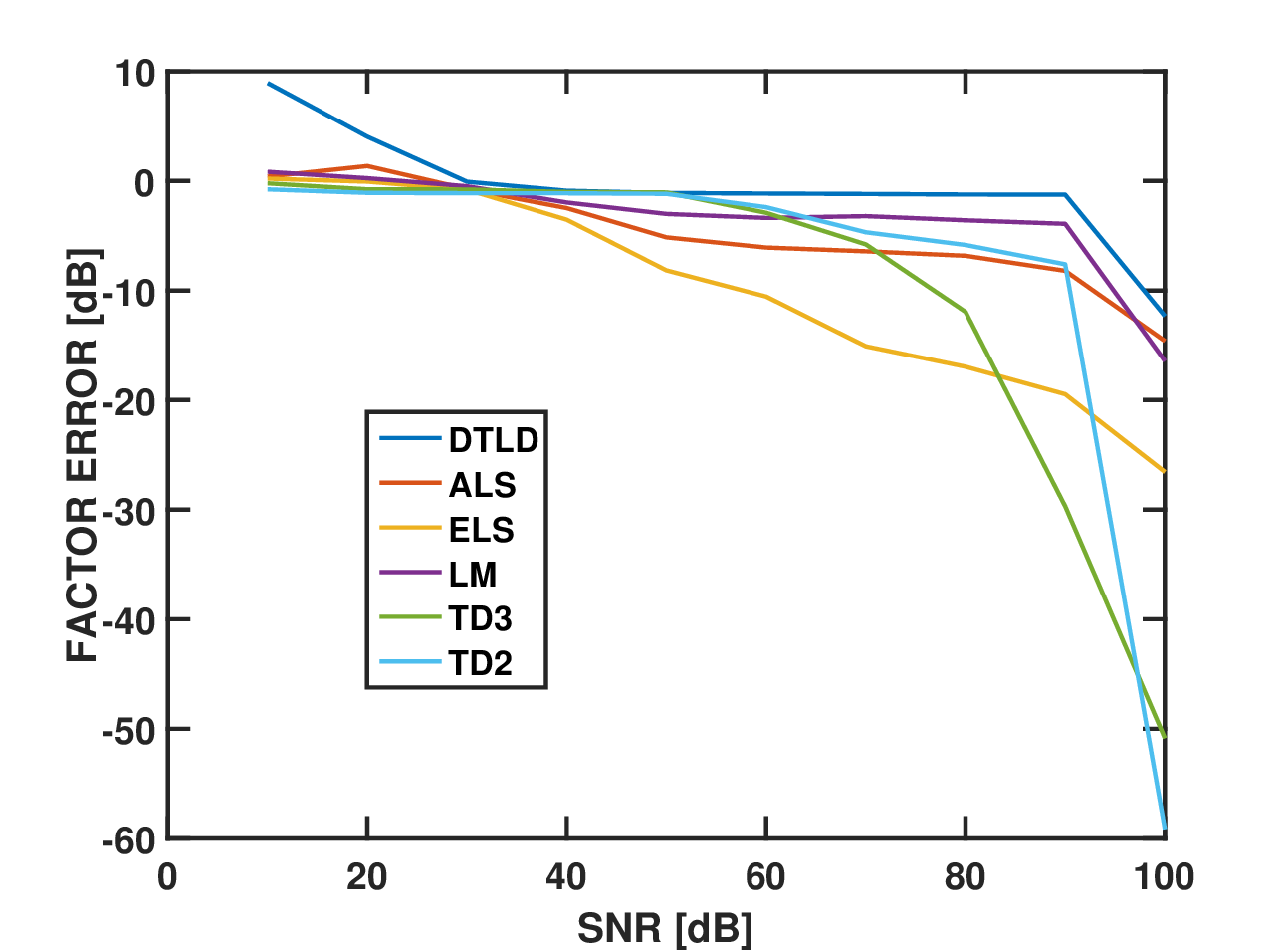}\label{f2} %{bcd_I12R3_medsae}
            \caption{Median fitting error and median factor estimation error of six CP decomposition algorithms
            as function of for as a function of the input SNR for $c=0.99$.}%
    %        (upper diagrams) and for SNR=30dB and varying $c$ (lower diagrams).}%\label{fig_MSAE}
\end{figure}

The decomposition was performed 100 times, each time with a new tensor and new additive noise added to the tensor.
As the measure of success, we considered two criteria: (1) the median fitting error between the noisy tensor and its CP decomposition model, and (2)
median error between the estimated and theoretical factor matrices. In the latter case, one must solve the permutation and scaling ambiguity of the columns of the estimated factor matrices to
match the theoretical ones. We computed median of sum of squared angular estimation errors between columns of the estimated and theoretical factor matrices.
The resulting criteria are plotted versus the input tensor signal-to-noise ratio (SNR) in Figure \ref{f2}.

Both criteria are, in general, decreasing functions of the SNR. While the fitting error converges quite monotonously, the factor error
remains high until SNR is as high as 80 dB.
According to the fitting error, the best method seems to be the ALS-ELS, if LM with large number of iterations is not counted. It is because the method more often than its competitors, converged close to the minimum mean square fitting error. On the other hand, we can see in the second diagram, that lower fitting error may not always imply lower error in the estimated factor matrices. For SNRs higher than 80 dB, TEDIA 2 and TEDIA 3 achieve lower factor error than the other methods. In this region of the SNR, the other methods are not that good, because DTLD is not
sufficiently good initialization for them.

\begin{figure}
            \centering
            \includegraphics[width=0.88\linewidth]{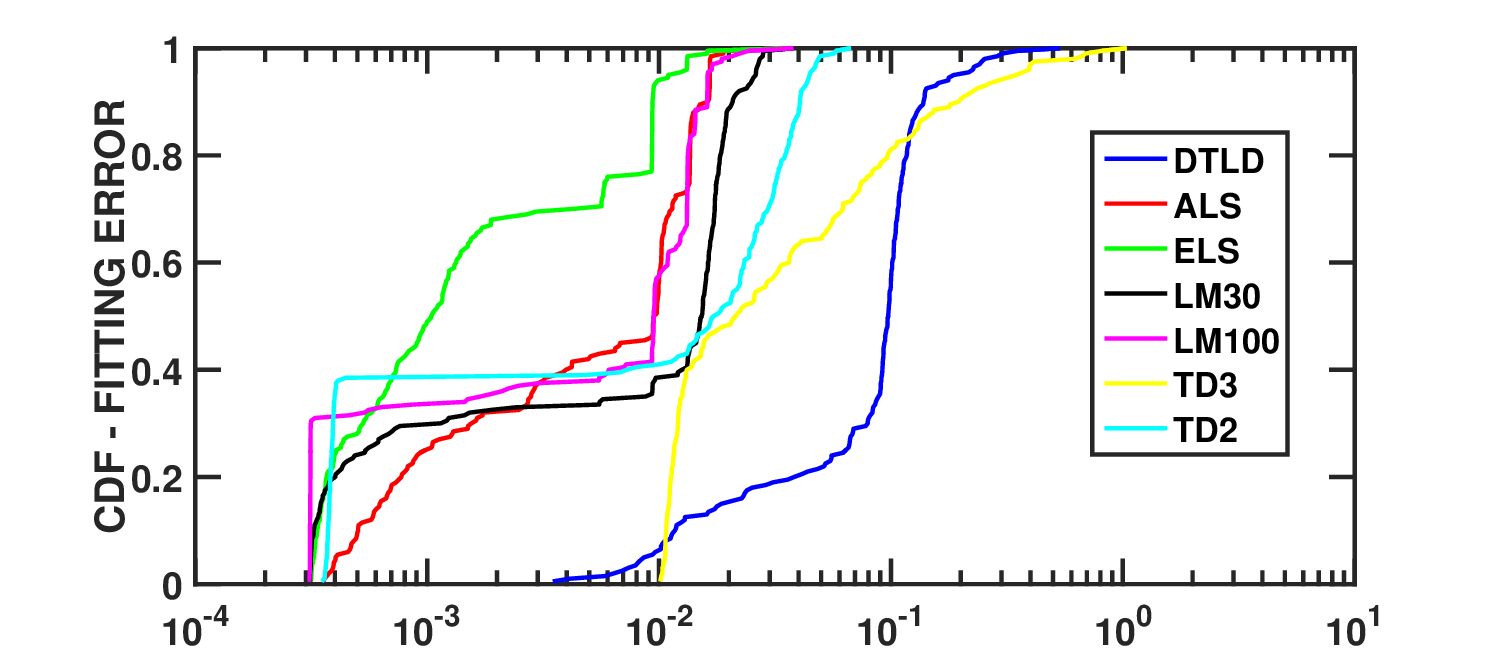}\\
            \includegraphics[width=0.88\linewidth]{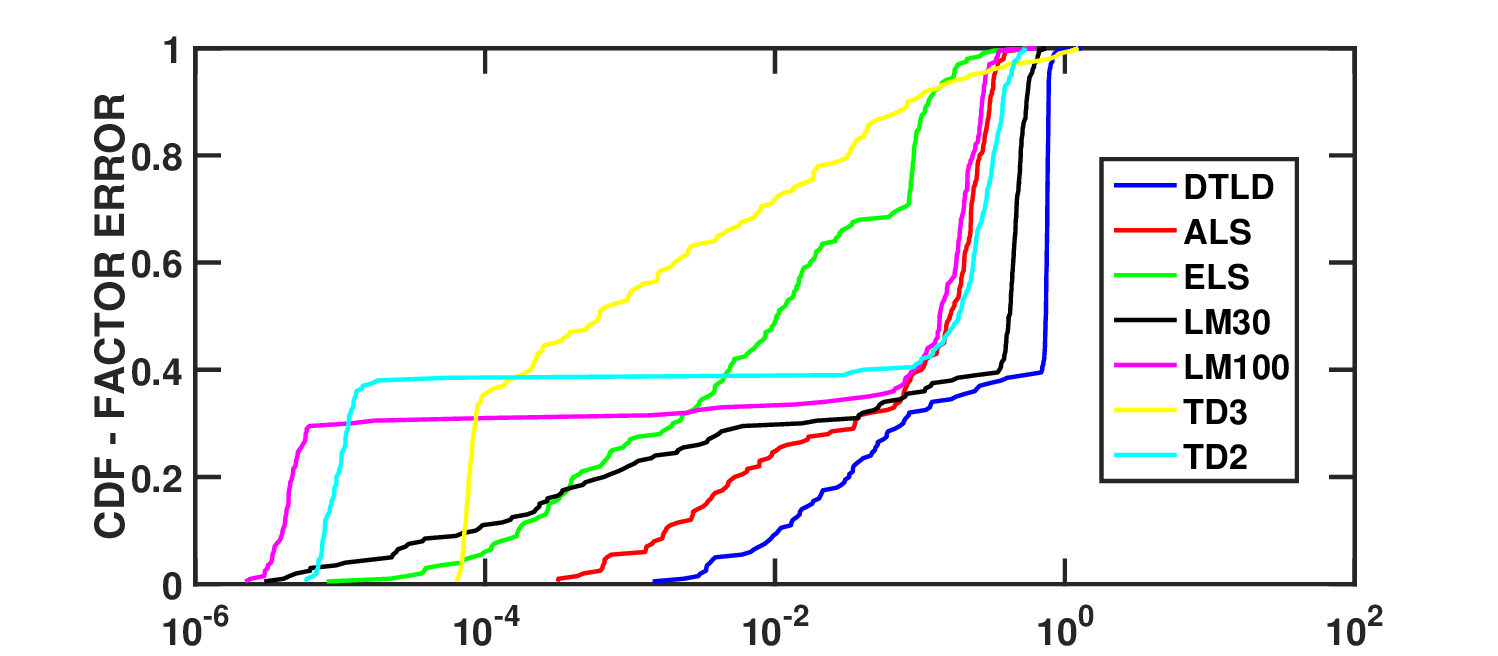}\vspace*{5mm} \\%{bcd_I12R3_medsae}
            \caption{Cumulative distribution function of the fitting error (upper diagram) and of the factor error (lower diagram) for the previous example with SNR=90 dB.}\label{f3}
    %        (upper diagrams) and for SNR=30dB and varying $c$ (lower diagrams).}%\label{fig_MSAE}
\end{figure}

To confirm the above observations, we considered the value SNR=90 dB, and studied cumulative distribution function of the fitting errors and the factor errors. Results are shown in Figure \ref{f3}. The LM method was studied twice, once with max. of 30 iterations, and second time with maximum of 100 iterations. With this number of iterations, the method was the best one in 30\% trials only. The median fitting error was the lowest for TEDIA 3, unless the number of iterations of LM is increased to cca. 1000 (but the method becomes slow).

We can conclude that TEDIA 2 and TEDIA 3 are good in decomposing tensors with a zero or a small noise.
In difficult scenarios, where the other method frequently fail due to their convergence to false minima of the cost function, TEDIA 2 and TEDIA 3 might do better.
In practice, of course, it is possible to combine the methods with the LM to achieve uniformly optimum results.

%In the second example, referred to as a double bottleneck, we keep the same $\bA$ and $\bB$ as in the first example, and $\bC=\bC_0$.

\subsection{Example 2: Approximate Joint Block AJD}

We have compared performance of seven approximate joint block AJD algorithms:
(1) U-WEDGE completed by clustering of rows of a demixing
matrix: this algorithm is blind to the assumed block structure.
This algorithm is used to initialize all subsequent
ones; (2) algorithm JBD NCG [6], %(3) the ad hoc algorithm
%of Ghennioui et al [4],
(3) the LLAJD algorithm \cite{Lahat},
%(5) FITJBD \cite{eusipco09},
and finally (4) the block two-sided TEDIA algorithm proposed in this paper.

We consider ten target matrices, each having four diagonal
blocks of the size $10\times 10$. The blocks were taken at random,
different at each simulation trial: each block is taken as the product
$\bX_{km}\bX_{km}^T$, where $\bX_{km}$ is Gaussian-distributed with zero mean and
variance one, mutually independent entries and independent in different slices.
Here, $k$ is an index of block, $k=1,2,3,4$ and $m$ is index of slice, $m=1,\ldots,10$.
Thus the resultant core tensor $\tS$ has dimension $40\times 40\times 10$
and is composed of four blocks of the size $10\times 10\times 10$.
The noisy tensors in the simulations are not obtained by adding a noise but they are built of
sample covariance matrices of $T$ random vectors having the required
theoretical covariances. To be specific, let $\bR_m$ be the $m-$th slice of the original tensor,
$\bR_m=\bX_{m}\bX_{m}^T$, $\bX_{m}$ being block diagonal with blocks $\bX_{km}$, then
the corresponding noisy tensor slice is
$\widehat\bR_m=\frac{1}{T} \bX_{m}\bY_T\bY_T^T\bX_{m}^T$, where $\bY_T$ is a random
matrix of the size $N\times T$ with Gaussian i.i.d. entries of zero mean and unit variance.
The resultant tensor is block dominant, but not exactly block diagonal.

The mixing matrix A was taken at random, also new
in each simulation trial. We compute it from a random unitary matrix $\bA_0$
as $\bA=c\bA_0(:,1){\bf 1}_{1\times 40}+\sqrt{1-c^2}\bA_0$, like in Section 7.1,
to obtain mixing matrix with collinear columns. We set $c=0.8$.
The mixture is the tensor $\tT=\tS \times_{\rm 1} \bA \times_{\rm 2}\bA$.
The block structure of the core tensor $\tS$ implies the tensor $\tT$ decomposition
as a sum
\begin{equation}
\tT=\tT_1+\tT_2+\tT_3+\tT_4~.\label{decomp}
\end{equation}
Each of the tensor $\tT_i$, $i=1,2,3,4$, has the size $40\times 40\times 10$ and multilinear rank $(10,10,10)$.
The block-term decomposition has several indeterminacies, e.g. the bases of the
independent subspaces can be quite arbitrary, but the decomposition (\ref{decomp})
is unique up to the order of the terms in the sum. Therefore, we shall measure success of the
approximate joint block diagonalization by mean square errors of appropriately sorted estimates of $\tT_i$,
$i=1,2,3,4$.

\begin{figure}
            \centering
            \includegraphics[width=0.78\linewidth]{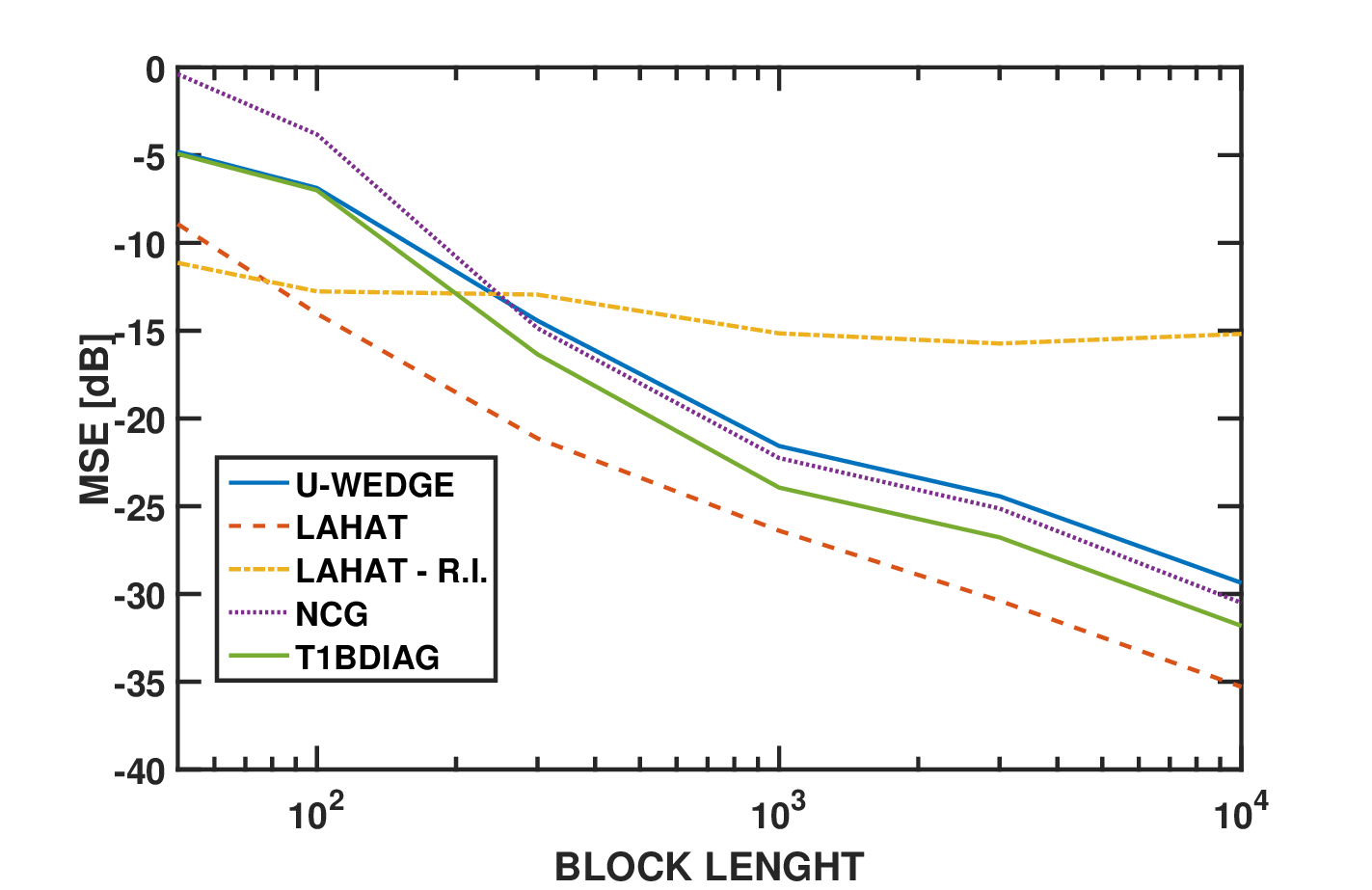}
            \caption{Mean square fitting error in dB of separated tensors
             of multilinear rank (10,10,10) for UWEDGE, JBD \cite{Lahat}, JBD \cite{Lahat} with random initialization,
             NCG algorithm \cite{Nion} and block TEDIA versus parameter $T$.}\label{JBD}
\end{figure}

We studied performance of four JBD algorithms: UWEDGE followed by collecting
the columns so that the block structure is revealed, JBD of \cite{Lahat} initialized by
the outcome of UWEDGE, JBD of Lahat et.al. \cite{Lahat} with default (random) initialization,
NGG algorithm of \cite{Nion}, and Block TEDIA initialized by the outcome of UWEDGE.
Results are presented in Figure \ref{JBD}.
First, we observe that the best performance is obtained by the algorithm of Lahat \cite{Lahat} that has been initialized by the outcome of UWEDGE. It is because the data generation model is in accord with this algorithm. The second best algorithm is the block TEDIA. The running times were 0.32 s, 0.99 s, 4.1 s and 0.98 s for UWEDGE, JBD \cite{Lahat}, NCG \cite{Nion} and TEDIA, respectively.

\subsection{Example 3: Three-sided block diagonalization} % / block term decomposition}

The initial tensor of the size $20\times 20\times 20$ is block
diagonal, with four random blocks along its main diagonal, each
of the size $5\times 5\times 5$. These blocks were computed as
a diagonal tensor having 1's on its main diagonal plus Gaussian random noise with zero mean and unit variance.

The initial tensor is the desired core tensor $\tE$.
The factor matrices $\bA$, $\bB$, $\bC$ were taken at random, as in the previous example.
We compute them from a random unitary matrices $\bA_0$, $\bB_0$, $\bC_0$
as $\bA= c\bA_0(:,1){\bf 1}_{1\times 20}+\sqrt{1-c^2}\bA_0$, e.t.c, like in Section 7.1,
We set $c=0.5$ and $c=0.9$, respectively.
The mixture is the tensor $\tT=\tE\times_1 \bA\times_2\bB \times_3 \bC$.
The block structure of the core tensor $\tS$ implies the tensor $\tT$ decomposition
as in (\ref{decomp}).

Each of the tensor $\tT_i$, $i=1,2,3,4$, has multilinear rank $(5,5,5)$.
The block-term decomposition has several indeterminacies, e.g. the bases of the
independent subspaces can be quite arbitrary, but the decomposition (\ref{decomp})
is unique up to the order of the terms in the sum. Therefore, we shall measure success of the
approximate joint block diagonalization by mean square errors of appropriately sorted estimates of $\tT_i$,
$i=1,2,3,4$.
A Gaussian noise is added to the tensor $\tT$ according to pre-specified
SNR values. %An additive noise to the observed tensor $\tT$.
%the noise tensor is taken at random and symmetrized with respect to the first two coordinates, and
%scaled to achieve desired signal-to-noise ratio (SNR), which
%is defined as 10 log10 of Frobenius norm of the mixtures divided
%by Frobenius norm of the noise.

We have tested four BTD algorithms: (1) Blind TEDIA with 1000 iterations, i.e. three-sided diagonalization followed by collecting
the columns so that the block structure is revealed, (2) Fixed block-size TEDIA with 5000 iterations and random initialization
(3) Fixed block-size TEDIA with 1000 iterations after being initialized by outcome of the blind TEDIA
(4) Block Alternating Least Squares with 100 iterations after it is initialized by the blind TEDIA.
Results of 100 independent trials are presented in Figure 6.

\begin{figure}
            \centering
            \includegraphics[width=0.48\linewidth]{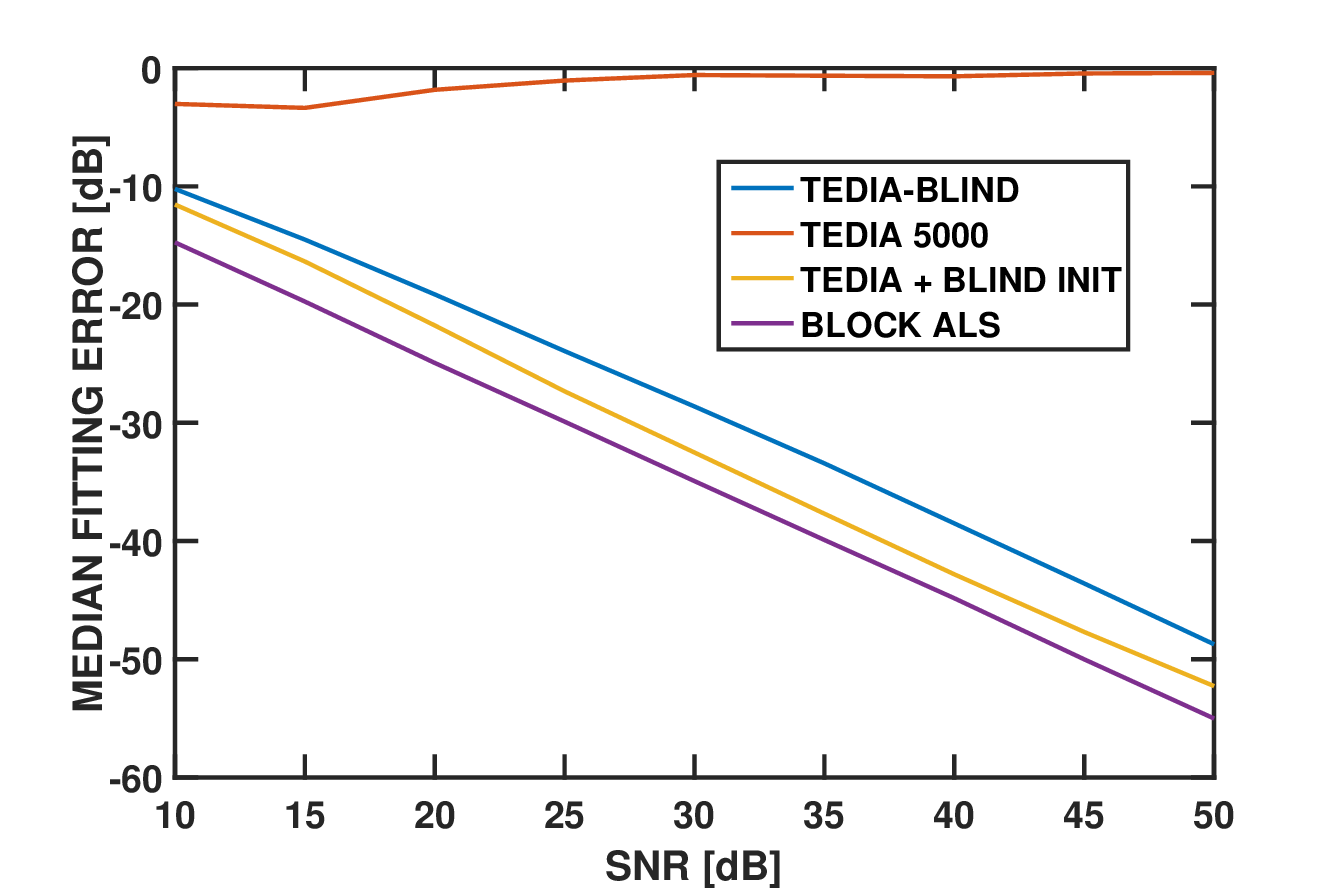}\quad%{bcd_I12R3_msae}
            \includegraphics[width=0.48\linewidth]{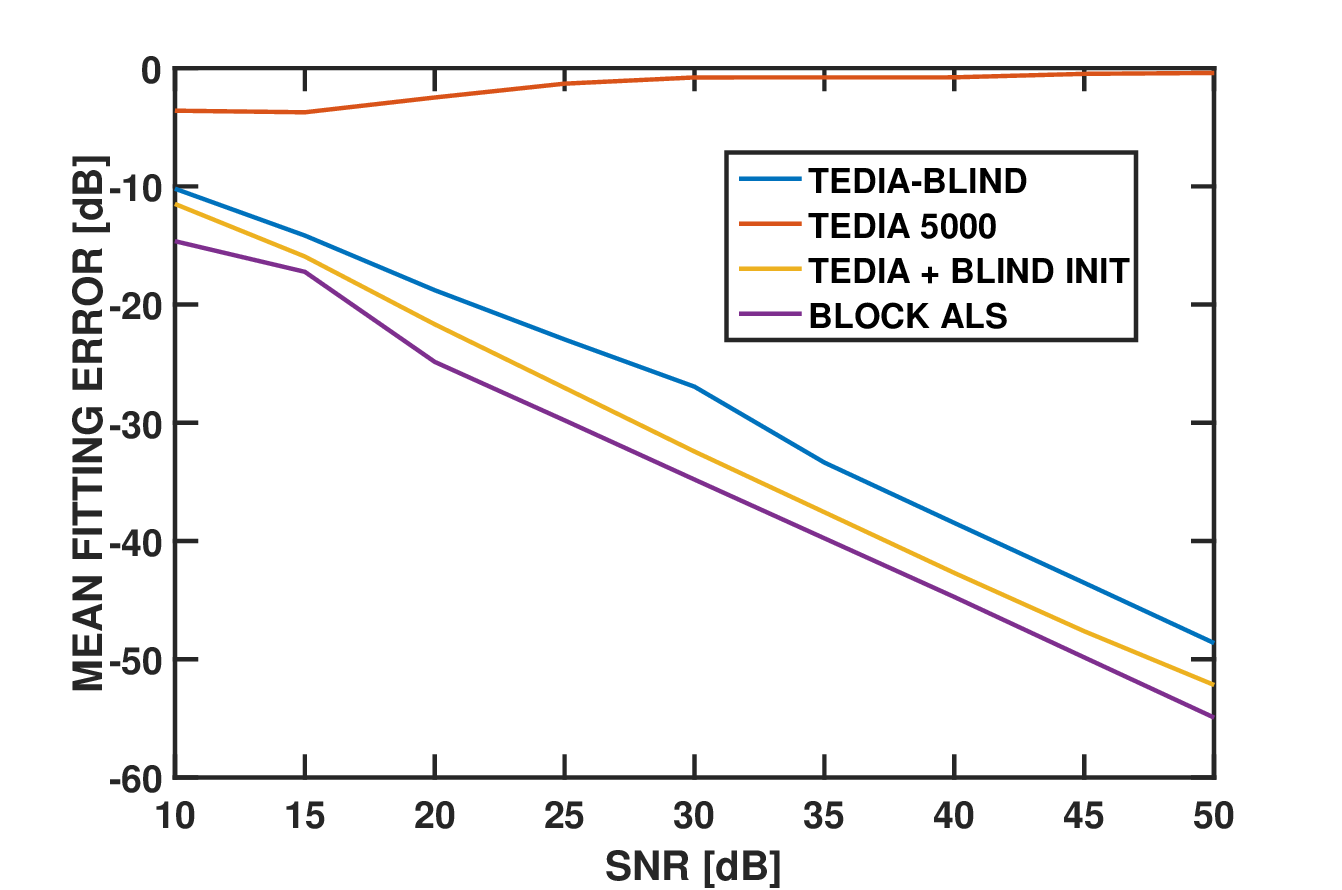}\vspace*{5mm} %{bcd_I12R3_medsae}
            \includegraphics[width=0.48\linewidth]{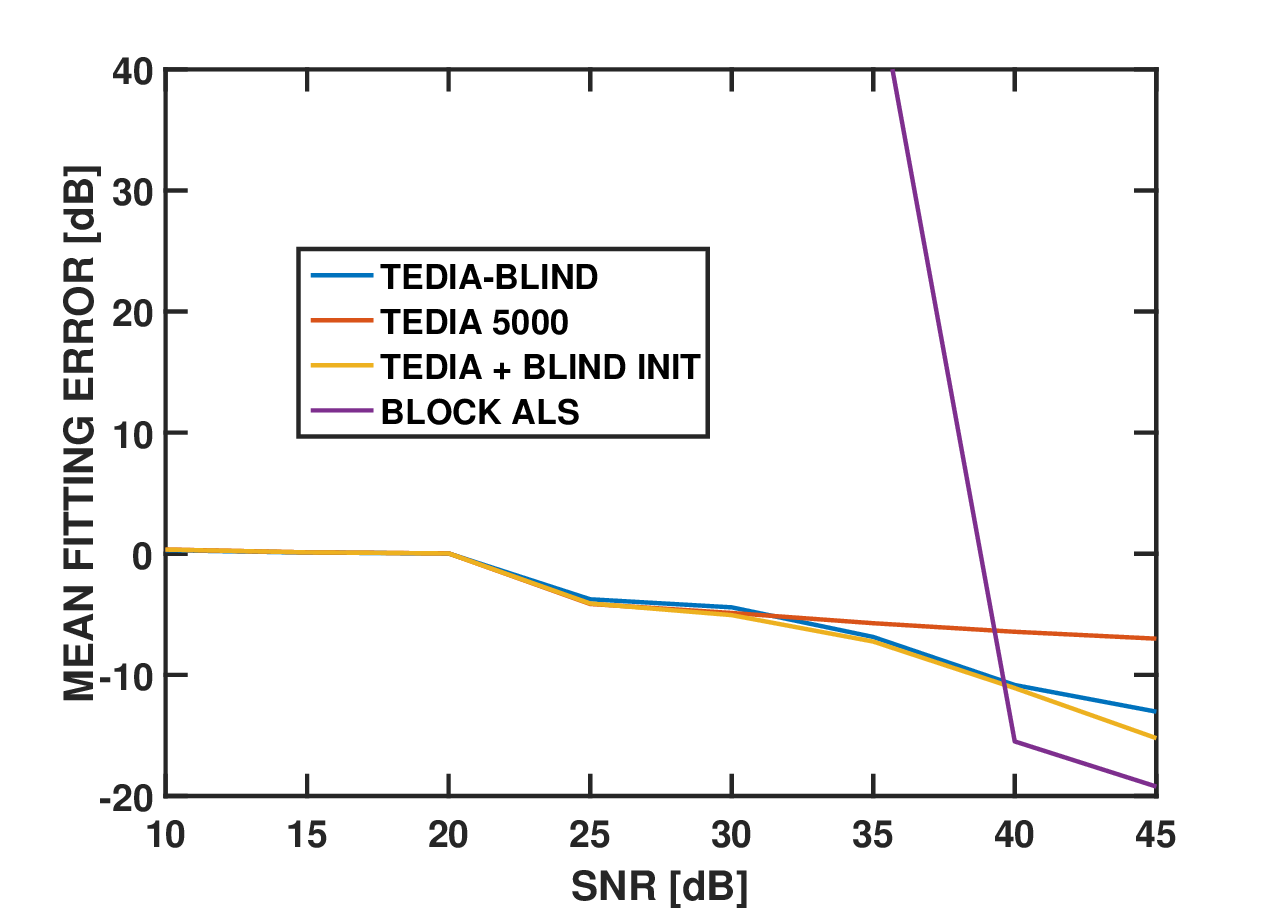}\quad%{bcd_I12R3_msae}
            \includegraphics[width=0.48\linewidth]{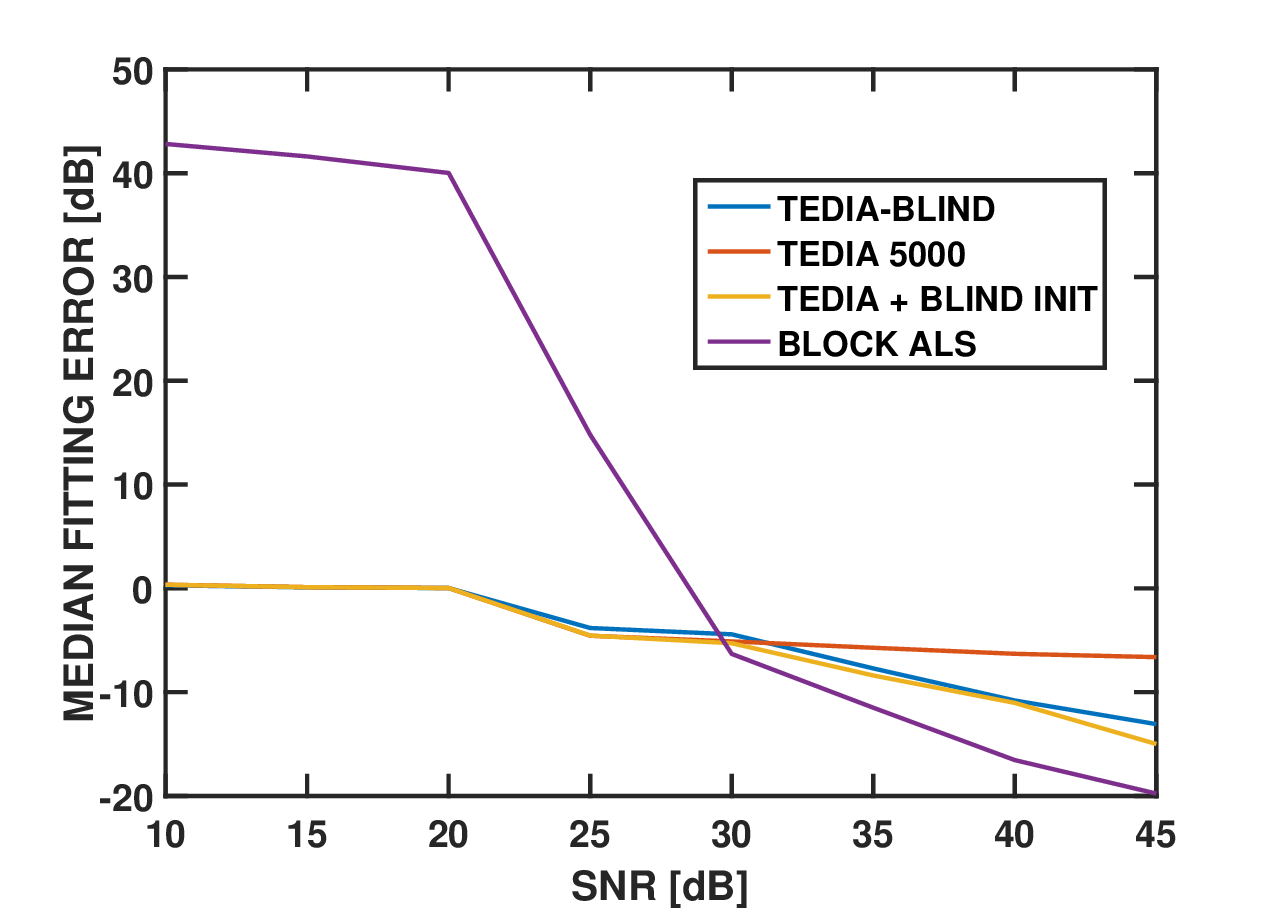}%{bcd_I12R3_medsae}
            \caption{Mean and median error in dB of separated tensors
             of multilinear rank (5,5,5) for (1) blind TEDIA with 1000 iterations, (2) block TEDIA with
             5000 iterations, (3) block TEDIA with 1000 iterations but initialized by outcome of the blind TEDIA, and
             (4) block ALS with 100 iterations initialized by the blind TEDIA. Upper diagram: $c=0.5$, lower digrams: $c=0.9$.}%\label{fig_MSAE}
\end{figure}
We can see that convergence of block TEDIA is relatively slow, 5000 iterations are not enough unless the algorithm is
initialized properly, e.g. by the outcome of the blind TEDIA. All three of these algorithms work relatively well even at low SNR's.
The behavior of the block ALS \cite{BTD3} is even more erratic in the difficult scenario with $c=0.9$, as we can see on the difference between the median and mean fitting error. The algorithm was initialized by the outcome of the blind TEDIA. The performance is good if the input SNR is sufficiently high:
40 dB. If the algorithm is initialized randomly, it usually fails.

Note that one run of the blind TEDIA (1000 iterations) takes 1.36 second, the additional 1000 iterations of the block TEDIA
requires additional 1.44 second. The 5000 iterations of the block TEDIA takes 9.35 seconds, and one run of the block ALS requires 71.9 seconds: it is very slow compared to TEDIA.

\section{Conclusions}

TEDIA is the technique of non-orthogonal tensor diagonalization
and block diagonalization. In difficult scenarios, it can outperform traditional methods of CP tensor
decomposition such as the alternating least squares (ALS), ALS with the enhanced line search, and Levenberg-Marquardt
method. The main reason for the success of TEDIA is that it does not suffer from many local minima of the cost function, unlike the traditional methods.
In the area of the block term decomposition, the situation is similar. We showed that TEDIA allows to fit the assumed block structure of the tensor directly, but sometimes it is useful to begin the separation with the blind TEDIA first.

Potential applications can be found in DS-CDMA systems
or in tensor deconvolution, in particular in feature extraction
and other areas.

Matlab code of the proposed technique is posted on the web page of
the first author.

\section*{Appendix A}

In this Appendix we derive the expression (\ref{gradient}) for $\bG_A$. The other gradients,
$\bG_B$ and $\bG_C$ follow from the symmetry of the problem.

Let
$$
\tE^\prime = \tE\times_1 \bA~.
$$
The $(k,\ell,m)$-th element of the tensor is
$$
e_{k\ell m}^\prime=\sum_{\alpha=1}^N
e_{\alpha\ell m} a_{\alpha k}~.
$$
Then
\begin{eqnarray*}
\frac{\partial e_{k\ell m}^\prime}{\partial a_{ij}}&=&
 \sum_{\alpha} e_{\alpha\ell m}
\delta_{\alpha i}\delta_{kj}= e_{i\ell m} \delta_{kj}
\end{eqnarray*}
and
\begin{eqnarray*}
\|\mbox{off}_3(\bE^\prime)\|^2 &=& \sum_{(k,\ell,m)\ne(k,k,k)}(e_{k\ell m}^\prime)^2\\
\frac{\partial\|\mbox{off}_3(\bE^\prime)\|^2}{\partial a_{ij}} &=& 2 \sum_{(k,\ell,m)\ne(k,k,k)}e_{k\ell m}\frac{\partial e^\prime_{k\ell m}}{\partial a_{ij}} \\ &=&\sum_{(k,\ell,m)\ne(k,k,k)}e_{k\ell m}e_{i\ell m} \delta_{kj}\\ &=&
\sum_{(\ell,m)\ne(j,j)}e_{j\ell m}e_{i\ell m}\\ &=&
\sum_{\ell,m=1}^N [\mbox{off}(\tE)]_{j\ell m}e_{i\ell m} ~.
\end{eqnarray*}

\end{document}